\documentclass[12pt]{amsart}

\usepackage{mathrsfs}

\usepackage{amsmath,amsthm,amssymb}
\font\teneufm=eufm10
\font\seveneufm=eufm7
\font\fiveeufm=eufm5
\newfam\frakturfam

\textfont\frakturfam=\teneufm
\scriptfont\frakturfam=\seveneufm
\scriptscriptfont\frakturfam=\fiveeufm

\newtheorem{lm}{Lemma}
\newtheorem{theor}{Theorem}
\newtheorem{co}{Corollary}

\def\bee{\begin{eqnarray}}
\def\bes{\begin{eqnarray*}}
\def\eee{\end{eqnarray}}
\def\ees{\end{eqnarray*}}

\def\a{\alpha}
\def\b{\beta}
\def\g{\gamma}

\def\Proof{{\sl Proof.}\ }

\title{A Non-linearizable Action of $S_3$ on Free Metabelian Lie and Associative Algebras}

\begin{document}
\date{}
\maketitle

\begin{center}

{\bf Altyngul Naurazbekova}\footnote{
L.N. Gumilyov Eurasian National University,
Astana, 010008, Kazakhstan,
e-mail: {\em naurazbekova82@gmail.com}}
and 
{\bf Ualbai Umirbaev}\footnote{Department of Mathematics,
Wayne State University,
Detroit, MI 48202, USA; Center for Theoretical and Applied Algebra, 
Al-Farabi Kazakh National University, Almaty, 050040, Kazakhstan,
e-mail: {\em umirbaev@wayne.edu}}

\end{center}

\begin{abstract}
Let $M$ be the free metabelian Lie algebra of rank $4$ over the field of complex numbers $\mathbb{C}$. We construct  a non-linearizable action of the symmetric group $S_3$, the smallest nonabelian group,   on $M$ by automorphisms. We also  construct  a non-linearizable action of  $S_3$ on free metabelian associative algebras of rank $4$. 
\end{abstract}

\


\tableofcontents

\section{Introduction}

\hspace*{\parindent}

The well-known Jung--van der Kulk Theorem \cite{Jung,Kulk} states that every automorphism
of the polynomial algebra $K[x,y]$ in two
variables over an arbitrary field $K$ is tame. Moreover, its automorphism group admits an amalgamated free product decomposition \cite{Kulk,Shafar}, namely
\bes
\mathrm{Aut}(K[x,y]) \cong \mathrm{Aff}_2(K) \ast_{H} \mathrm{Tr}_2,
\ees
where $\mathrm{Aff}_2(K)$ denotes the subgroup of affine automorphisms, $\mathrm{Tr}_2$ denotes the subgroup of triangular automorphisms, and $H=\mathrm{Aff}_2(K)\cap \mathrm{Tr}_2$.

In 1979, T. Kambayashi \cite{Kam}, using the amalgamated product structure of $\mathrm{Aut}(K[x,y])$, proved that any algebraic subgroup of $\mathrm{Aut}(K[x,y])$ is conjugate to a subgroup of affine automorphisms or triangular automorphisms. In particular, in the case of characteristic zero any action of a reductive group on $\mathbb{A}^2$ is linearizable, i.e., conjugate to a linear action by a polynomial automorphism of $\mathbb{A}^2$. In 1974, M. Suzuki proved \cite{Su} the same result for the complex affine space $\mathbb{A}^2$ using geometric methods. 

An analogue of this question for $\mathbb{A}^n$, $n\geq 3$, is well known as Kambayashi's Linearization Problem \cite{Essen}.  In the case $n=3$, the answer remains affirmative for certain special cases \cite{KKMR,Popov01}. The first examples of nonlinearizable actions  were given by G.W. Schwarz \cite{Schwarz89} and this work initiated an exciting development. Examples of non-linearizable actions of finite groups appeared in \cite{MMP98,MP95,Med}. Finally, G. Freudenburg and L. Moser-Jauslin constructed a non-linearizable action of $S_3$ on $\mathbb{A}^4$ \cite{FM02}. 

Nevertheless, an important question concerning the linearizability of finite-order polynomial automorphisms of the complex affine space $\mathbb{A}^n$ remains open for all $n\geq 3$ \cite{KS}. 

It is well known that the automorphism groups of  free associative algebras in two variables over a field \cite{Czer,Makar}, as well as that of  free Poisson algebras in two variables over a field of characteristic zero \cite{MLTU}, are isomorphic to the automorphism group of the polynomial algebras in two variables. In 2002, P. Cohn \cite{Cohn02} proved that every finite group of automorphisms of a free associative algebra in two variables over a field of characteristic zero is linearizable. 
In 2019, A. Alimbaev, A. Naurazbekova, and D. Kozybaev \cite{ANZ} showed that any algebraic subgroup of the group of automorphisms of free associative algebras, free Poisson algebras, and free left-symmetric algebras in two variables is conjugate to a subgroup of affine automorphisms or triangular automorphisms. Consequently, in the case of characteristic zero  any action of a reductive group on any of these algebras by automorphisms is linearizable. 

In 1964, P. Cohn \cite{Cohn1} proved that all automorphisms of finitely generated free Lie algebras are tame. Automorphisms of free Lie algebras in two variables are linear. M.A. Shevelin proved \cite{Shevelin11,Shevelin12} that any finite group of automorphisms of free Lie algebras in three variables over a field of characteristic zero is linearizable.  In 2020, it was proved \cite{ANU20} that the  automorphism group of a free Lie algebra $L_3=\mathrm{Lie}\langle x,y,z\rangle$ in three variables $x,y,z$ admits an amalgamated free product structure
\bes
\mathrm{Aut}(L_3)=\mathrm{GL}_3(K) \ast_H T_3,
\ees
where $\mathrm{GL}_3(K)$ is the subgroup of all linear automorphisms, $T_3$ is the subgroup of all automorphisms of the form
\bes
x\mapsto \a x+f, \ y\mapsto  \b_2 y+\b_3 z, \  z\mapsto \g_2 y+\g_3 z,
\ees
where $\a,\b_i,\g_i\in K$, $\a(\b_2\g_3-\b_3\g_2)\neq 0$, $f\in \mathrm{Lie}\langle y,z\rangle$, and $H=\mathrm{GL}_3(K)\cap T_3$. This result allows to prove that any algebraic subgroup of the automorphism group  of a free Lie algebra in three variables is conjugate to a subgroup of $\mathrm{GL}_3(K)$ or $T_3$. Consequently, any action of a reductive group on a free Lie algebra in three variables over a field of characteristic zero by automorphisms is linearizable.

In 2014, M.A. Shevelin \cite{Shevelin14} constructed an interesting example of a subgroup of $\mathrm{Aut}(L_n)$, $n\geq 3$, over the field of complex numbers $\mathbb{C}$, which is isomorphic to a subgroup of $\mathbb{C}^*$ but is not conjugate to a subgroup of $\mathrm{GL}_n(\mathbb{C})$. Modifying this example, V.G. Bardakov and M.V. Neshchadim \cite{BN} proved that for any prime $p$ there exists a continuum of pairwise non-conjugate subgroups of the automorphism group  $\mathrm{Aut}(\mathbb{C}[x,y])$ of the polynomial algebra $\mathbb{C}[x,y]$ such that each of them is isomorphic to the quasi-cyclic $p$-group
$$\mathbb{C}_{p^\infty}=\{x\in \mathbb{C} | x^{p^n}=1 \  \mathrm{for} \  \mathrm{some} \  n\in \mathbb{N}\}$$
and each of them is not conjugate in $\mathrm{Aut}(\mathbb{C}[x,y])$ to a linear subgroup of $\mathrm{Aut}(\mathbb{C}[x,y])$.

Let $M_n$ be the free metabelian Lie algebra  freely generated by  $x_1,x_2,\ldots,x_n$ over an arbitrary field $K$. It is well known that every nontrivial exponential automorphism of $M_2$ is wild \cite{Art2,Shm}. In 1992, V. Drensky \cite{Dr} proved that the exponential automorphism $\mathrm{exp}(\mathrm{ad}[x_1,x_2])$ of $M_3$ is wild. 
R. Nauryzbaev \cite{Naur} proved this automorphism of $M_3$ is not even almost tame. It was shown in \cite{NDU} that both the group of all tame automorphisms and the group of all almost tame automorphisms of  $M_3$ admit structures of amalgamated free products.
In 2008, V. Roman'kov \cite{Rom} proved that the group of automorphisms of $M_3$ cannot be finitely generated modulo the subgroup generated by all exponential and almost tame automorphisms.  

In 1993, Bryant and Drensky proved \cite{BD93-2} that every automorphism of $M_n$ of rank $n\geq 4$ over a field $K$ of characteristic zero is {\em approximately tame} ( see the definitions in \cite{SU25}); that is, it can be approximated by tame automorphisms with respect to power series topology \cite{BD93-2,SU25}. In contrast, all wild automorphisms of $M_3$ are {\em absolutely wild} \cite{SU25}, meaning that they cannot be approximated by tame automorphisms.

Recently, U. Umirbaev \cite{UUU6,UUU7} proved that if $n\geq 4$ then every automorphism of $M_n$ is almost tame. Moreover, he showed that, for $n\geq 4$, the group  $\mathrm{Aut}(M_n)$ of all automorphisms of the free metabelian Lie algebra $M_n$ over a field $K$ of characteristic not equal to $3$ is generated by its linear automorphisms, together with a single quadratic and a single cubic automorphism.

In this paper, we construct a non-linearizable action of $S_3$  on the free metabelian Lie algebra $M$ of rank $4$ over the field of complex numbers $\mathbb{C}$ by automorphisms. We also  construct a similar action for free metabelian associative algebras. 

The question of whether finite groups of automorphisms of free associative algebras and free Lie algebras are linearizable remains open. Recall that finite linear groups of automorphisms of free associative algebras were studied by V. Kharchenko \cite{Khar}.

The paper is organized as follows. In Section 2, we recall the definition and some properties of the non-linearizable action of $S_3$ on the polynomial algebra from \cite{FM02,FM04}.  In Section 3, we review some results on the structure of free metabelian Lie algebras and their Fox derivatives. 
In Section 4, following some ideas from \cite{FM02,FM04}, we construct an action of  $S_3$ on the free metabelian Lie  algebra  $M$. In Section 5, we prove that this action of $S_3$ is not linearizable. In Section 6, we construct a similar action of  $S_3$ on the free metabelian associative algebra  $MA$ and prove that it is also not linearizable.

\section{The Action of $S_3$ on Polynomial Algebras}
\hspace*{\parindent}

Let $\mathbb{C}[a,b,x,y]$ be the polynomial algebra over $\mathbb{C}$ in the variables $a,b,x,y$. We recall the non-linearizable action of $S_3$ on $\mathbb{C}[a,b,x,y]$ constructed in \cite{FM02}.

The group $S_3$ is generated by $\sigma$ and $\tau$, subject to the relations
\bee\label{f1}
\sigma^3=1,  \ \tau^2=1, \  \text{and}  \ \tau\sigma=\sigma^2\tau.
\eee

Denote by
\bes
\phi=(f,g,h,m)
\ees
or by
\bes
\phi=(
 \left[\begin{array}{cc}
 f  \\
 g  \\
\end{array}\right],
\left[\begin{array}{cc}
 h  \\
 m  \\
\end{array}\right])
\ees
the endomorphism $\phi$ of $\mathbb{C}[a,b,x,y]$ defined by $a\mapsto f, b\mapsto g, x\mapsto h, y\mapsto m$.

 Set
\bes
\sigma=(wa,w^{-1}b,x,y), \ \text{where} \ w=e^{2\pi i/3},
\ees
 and
\bes
\tau=(
 \left[\begin{array}{cc}
 b  \\
 a  \\
\end{array}\right],
N\left[\begin{array}{cc}
 x  \\
 y  \\
\end{array}\right]), \ \text{where} \ N=\left[\begin{array}{cc}
 -b^3  &  1+ab+(ab)^2  \\
 1-ab  &  a^3  \\
\end{array}\right]. 
\ees
Observe that 
\bes
J(\tau)=\begin{bmatrix}
\begin{matrix}
0 & 1\\
1 & 0
\end{matrix} & \begin{matrix}
		0 & 0 \\
		0 & 0
	\end{matrix} \\
	\begin{matrix}
		* & * \\
		* & *
	\end{matrix} & N
\end{bmatrix}.
\ees
The automorphisms $\sigma$ and $\tau$ satisfy the relations (\ref{f1}) and this determines an action of $S_3$ on $\mathbb{C}[a,b,x,y]$. Notice that $\sigma$ is linear and $\tau$ is linearizable \cite{MJ93,MMP96}, but the action of $S_3$ is not linearizable \cite{FM02}.

\section{Free metabelian Lie algebras}

\hspace*{\parindent}

Recall that a Lie algebra $L$ over a field $K$ is  called {\em abelian} if $[L,L]=0$ and {\em metabelian} if $[[L,L],[L,L]]=0$.

Let $M_n$ be the free metabelian Lie algebra over a field $K$ freely generated by  $x_1,x_2,\ldots,x_n$. It is easy to verify that the identity $[x,[y,z]]=[y,[x,z]]$
 holds for all $x,y\in M_n$ and  $z\in[M_n,M_n]$ . The left normed monomials
\bes
[x_{i_m},\ldots,[x_{i_2},[x_{i_1},x_i]]\ldots], \ i>i_1\leq i_2\leq\ldots \leq i_m, \ m\geq 0, 
\ees
form a linear basis of $M_n$ \cite{Bakh}.

Let $A_n=M_n/[M_n,M_n]$. The elements $y_i=x_i+[M_n,M_n]$, $1\leq i\leq n$, form a linear basis of $A_n$. The universal enveloping algebra $U(A_n)$ of $A_n$ is the polynomial algebra $K[y_1,y_2,\ldots,y_n]$ in $n$ variables over $K$. For any $v=\sum \alpha_{i_1,\ldots,i_n}y_1^{i_1}\ldots y_n^{i_n}$, $\alpha_{i_1,\ldots,i_n}\in K$,  and any $g\in [M_n,M_n]$ we put
\bes
vg=\sum_{i_1,i_2,\ldots,i_n} \alpha_{i_1,\ldots,i_n} (\mathrm{ad}_{x_1})^{i_1}\ldots (\mathrm{ad}_{x_n})^{i_n} (g),
\ees
where $\mathrm{ad}_{x_i} (t)=[x_i,t]$, $t\in M_n$. The space $[M_n,M_n]$ becomes a left $U(A_n)$-module with respect to this action.

Denote by $\hat{f}$ the homomorphic image of the element $f\in M_n$ under the natural homomorphism $M_n\rightarrow A_n\subset U(A_n)$. According to 
\cite{Shp, UUU1}, the Fox derivatives
$$\frac{\partial}{\partial x_j}: M_n\rightarrow U(A_n),\  1\leq j\leq n,$$
are uniquely defined as linear maps satisfying the conditions
$$\frac{\partial x_i}{\partial x_j}=\delta_{ij},\  \frac{\partial[f,g]}{\partial x_j}=\hat{f} \frac{\partial g}{\partial x_j}-\hat{g} \frac{\partial f}{\partial x_j},$$
where $\delta_{ij}$ is the Kronecker delta.

Let $\varphi=(f_1,\ldots,f_n)$ be an arbitrary $n$-tuple of elements of $M_n$. Usually $\varphi$ considered as an endomorphism of $M_n$ defined by $\varphi(x_i)=f_i$ for all $i$. The Jacobian matrix of this endomorphism is 
\bes
J(\varphi)=\left[\frac{\partial f_i}{\partial x_j}\right]_{1\leq i,j\leq n}=\left[\begin{array}{cc}
 \partial(f_1) \\
 \partial(f_2) \\
\ldots \\
\partial(f_n) 
\end{array}\right],   
\ees
where $\partial(f_i)=(\frac{\partial f_i}{\partial x_1},\frac{\partial f_i}{\partial x_2}, \ldots, \frac{\partial f_i}{\partial x_n})$ for all $1\leq i\leq n$.   Recall that $\varphi$ is an automorphism if and only if $J(\varphi)$ is invertible over $U(A_n)$ \cite{Shp, UUU1}. 

\section{An Action of $S_3$ on Metabelian Lie Algebras}

\hspace*{\parindent}

Let $M=M\left\langle a,b,x,y\right\rangle$ be the free metabelian Lie algebra over the field of complex numbers $\mathbb{C}$ with free generators $a,b,x,y$.  We want to construct an action of $S_3$ on $M$ "similar" to the action of $S_3$ on  $\mathbb{C}[a,b,x,y]$. Similarity will be measured by their Jacobian matrices. Set 
\bes
\sigma=(wa,w^{-1}b,x,y), \ \text{where} \ w=e^{2\pi i/3},
\ees 
as in the case of polynomial algebras. Then $\sigma^3=1$. 

Denote $\mathrm{ad}_{a}$ and $\mathrm{ad}_{b}$ by $A$ and $B$, respectively. Let $C$ be the associative algebra generated by $A$ and $B$. This algebra is not commutative, but it is "almost commutative" in the following sense. Let $C_0$ be the ideal of $C$ generated by $A$ and $B$. Then $[A,B]C_0=0$. The elements
\bes
A^iB^j,\quad A^iB^j[A,B], \qquad i,j\geq 0,
\ees
form a linear basis of $C$. In fact, $C_0$ is a free algebra in the variety of associative algebras satisfying the identity $[x,y]z=0$. 

Let 
\bes
\tau=(
\left[\begin{array}{cc}
	b  \\
	a  \\
\end{array}\right],
S
\left[\begin{array}{cc}
	x  \\
	y  \\
\end{array}\right]),
\ees
 where
 \bes
 S=
 \begin{bmatrix}
 	-B^3+U[B,A] & 1+AB+A^2B^2+V[B,A]\\
 	1-AB+W[B,A] & A^3+Z[B,A]
 \end{bmatrix}, 
 \ees
$U,V,W,Z$ are arbitrary elements of $C$. Observe that 
\bes
J(\tau)=\begin{bmatrix}
	\begin{matrix}
		0 & 1\\
		1 & 0
	\end{matrix} & \begin{matrix}
		0 & 0 \\
		0 & 0
	\end{matrix} \\
	\begin{matrix}
		* & * \\
		* & *
	\end{matrix} & N
\end{bmatrix}
\ees
has the same form as $J(\tau)$ in the polynomial algebra. Thus $J(\tau)$ is invertible and $\tau$ is an automorphism for any choice of $U,V,W,Z\in C$. 

For any $F=F(A,B)\in C$ set $F^\ast = F(B,A)$. Note that the map $*$ defines an automorphism of $C$.

\begin{lm}\label{lu1} 
If 
\bes
V=1+AB \ \text{and} \ U=W=Z=0,
\ees
then $\tau^2=1$ and $\tau\sigma=\sigma^2\tau$.
\end{lm}
\Proof 
We have 
\bes
\tau^2=
(\begin{bmatrix}a\\ b\end{bmatrix}, \tau(S\begin{bmatrix}x\\y\end{bmatrix}))
=(\begin{bmatrix}a\\ b\end{bmatrix},S^\ast S\begin{bmatrix}x\\y\end{bmatrix}),
\ees
where
\bes
S^\ast S= \begin{bmatrix}
	-A^3-U^*[B,A] & 1+BA+B^2A^2-V^*[B,A]\\
	1-BA-W^*[B,A] & B^3-Z^*[B,A]
\end{bmatrix}\\
\cdot\begin{bmatrix}
	-B^3+U[B,A] & 1+AB+A^2B^2+V[B,A]\\
	1-AB+W[B,A] & A^3+Z[B,A]
\end{bmatrix}
=\begin{bmatrix}
	1+F_{11}[B,A] &   F_{12} [B,A] \\
	F_{21}[B,A] & 1+ F_{22}[B,A]
\end{bmatrix}
\ees
with
\bes
F_{11}=1+AB-A^3U+(1+AB+A^2B^2)W-V^\ast,\\
F_{12}=A^3+A^4B-A^3V+(1+AB+A^2B^2)Z-U^\ast,\\
F_{21}=(1-AB)U+B^3W-Z^\ast,\\
F_{22}=-1+A^2B^2+(1-AB)V+B^3Z-W^\ast.
\ees
Since $F_{ij}$, $1\leq i,j\leq 2$, appear to the left of $[B,A]$, the last four equalities can be considered in the polynomial algebra $\mathbb{C}[A,B]$.
It is easy to check that if  $V=1+AB$ and $U=W=Z=0$, then $F_{ij}=0$ for all $1\leq i,j\leq 2$. Consequently, $\tau^2=1$.

We also have 
\bes
\tau\sigma=(\begin{bmatrix}wb\\w^{-1}a\end{bmatrix},
\begin{bmatrix}
	-B^3 & 1+AB+A^2B^2\\
	1-AB & A^3
\end{bmatrix}
\begin{bmatrix}x\\y\end{bmatrix}+ \begin{bmatrix}
	U & V\\
	W & Z
\end{bmatrix}[B,A]
\begin{bmatrix}x\\ y \end{bmatrix})
\ees
and 
\bes
\sigma^2\tau=(\begin{bmatrix}wb\\w^{-1}a\end{bmatrix},
\begin{bmatrix}
	-B^3 & 1+AB+A^2B^2\\
	1-AB & A^3
\end{bmatrix}
\begin{bmatrix}x\\y\end{bmatrix}
+\b(\begin{bmatrix}
	U & V\\
	W & Z
\end{bmatrix})[B,A]
\begin{bmatrix} x\\ y \end{bmatrix}), 
\ees
where $\b$ is the automorphism of $C$ induced by $\sigma^2=(w^{-1}a, wb, x, y)$.

If  $V=1+AB$ and $U=W=Z=0$, then 
$$\begin{bmatrix} 	
U & V\\ 	W & Z 
\end{bmatrix}=
\b(\begin{bmatrix}
	U & V\\
	W & Z
\end{bmatrix}).$$ Consequently, $\tau\sigma=\sigma^2\tau$. $\Box$

Observe that, under the assumptions of Lemma \ref{lu1}, we have
\bee\label{fu1}
\tau=(b,a,-B^3x+(1+BA+B^2A^2)y, (1-AB)x+A^3y). 
\eee

From now until the end of Section 5, the symbol $\tau$ will denote exclusively this automorphism of the free metabelian Lie algebra $M$.    

\begin{co}\label{cu1}
The automorphisms $\sigma$ and $\tau$ of $M$ satisfy the relations (\ref{f1}). 
\end{co}

The statement of this corollary can also be verified directly, independently of Lemma \ref{lu1}. Corollary \ref{cu1} defines an action of $S_3$ on $M$. 

\section{Non-linearizability in the Case of Metabelian Lie Algebras}

\hspace*{\parindent}

In this section, we prove that the action of $S_3$ on $M$ defined above is not linearizable. 

An element $F$ of $M$ is called an $S_3$-invariant if $\sigma(F)=\tau(F)=F$, and an $S_3$-semi-invariant if $\sigma(F)=F$ and $\tau(F)=-F$.

Let $B_i $ be the linear span of all basis elements of degree $i$ in $x$ and $y$. This defines a grading 
$$M=B_0\oplus B_1\oplus \ldots \oplus B_m\oplus \ldots.$$
Since the action of $S_3$ preserves this grading, an element $h=\sum_{i=0}^{m}h_i\in M$, where $h_i\in B_i$, is $S_3$-invariant (respectively, $S_3$-semi-invariant) if and only if each homogeneous component $h_i$ is $S_3$-invariant (respectively, $S_3$-semi-invariant).

\begin{lm} \label{lm1}
 Let $h(a,b) \in \mathbb{C}[a,b]$. If $h(a,b)=h(wa,w^{-1}b)$, then $h(a,b)\in \mathbb{C}[a^3,b^3,ab]$.
\end{lm}
\Proof Let $h(a,b)=\sum_{m,n\geq 0} \gamma_{m,n}a^m b^n$, where $\gamma_{m,n} \in \mathbb{C}$ for all $m,n$. Since $h(a,b)=h(wa,w^{-1}b)$ it follows that
$$ \sum_{m,n\geq 0} \gamma_{m,n}a^m b^n=\sum_{m,n\geq 0} \gamma_{m,n}w^{m-n}a^m b^n.$$
Hence $w^{m-n}=1$ whenever $\gamma_{m,n}\neq 0$. Since $w$ is a primitive third root of unity,  $w^{m-n}=1$ if and only if $3\mid(m-n)$. If $m\geq n$, then $m-n=3k$ and
$$a^m b^n=a^{m-n}(a^nb^n)=(a^3)^k(ab)^n\in \mathbb{C}[a^3,b^3,ab].$$
If $n> m$, then $n-m=3k$ and
$$ a^mb^n=b^{n-m}(a^mb^m)=(b^3)^k(ab)^m\in \mathbb{C}[a^3,b^3,ab].$$ 
Consequently, $h(a,b)\in \mathbb{C}[a^3,b^3,ab]$. $\Box$

Denote by $\left\langle a-b\right\rangle$ the ideal of $M$ generated by $a-b$.

\begin{lm} \label{lm2} If $F\in B_1$ is an $S_3$-invariant, then 
\bes
 F=r((1-A)x+(1-A+A^2)y)\pmod{\left\langle a-b\right\rangle}, \ \text{where}  \ r=r(A) \notin \mathbb{C}^*.
\ees
\end{lm}
\Proof Let $F=F_x(A,B)x+F_y(A,B)y$, where $F_x(A,B),F_y(A,B)\in C$. It is well known that there exists a natural surjective homomorphism $\pi: C\rightarrow \mathbb{C}[a,b]$ uniquely determined by $\pi(A)=a$ and $\pi(B)=b$. For $h(A,B)\in C$  denote its image under $\pi$ by $h(a,b)$ . It is clear that
\bes
\frac{\partial F}{\partial z}=F_z(a,b), \ z\in\left\{x,y\right\}.
\ees

The invariance of $F$ under $\sigma$ implies
$$F=\sigma(F)=F_x(wA,w^{-1}B)x+F_y(wA,w^{-1}B)y.$$
Hence
\bes
F_z(a,b)=F_z(wa,w^{-1}b), \ z\in\left\{x,y\right\}.
\ees
By Lemma \ref{lm1}, it follows that $F_z(a,b) \in \mathbb{C}\left[a^3,b^3, ab\right]$, $z\in\left\{x,y\right\}$. 

The invariance  of $F$ under $\tau$  implies
\bes
F=\tau(F)=F_x(B,A)(-B^3x+(1+BA+B^2A^2)y)+F_y(B,A)((1-AB)x+A^3y). 
\ees
Hence
\bes
F_x(a,b)=F_x(b,a)(-b^3)+F_y(b,a)(1-ab),\\
F_y(a,b)=F_x(b,a)(1+ab+a^2b^2)+F_y(b,a)(a^3).
\ees
In the case $a=b$, the last two equalities imply that
\bes
F_x(a,a)(1+a^3)=F_y(a,a)(1-a^2),\\
F_x(a,a)(1+a^2+a^4)=F_y(a,a)(1-a^3),
\ees
where $F_z(a,a) \in \mathbb{C}\left[a^2,a^3\right]$, $z\in\left\{x,y\right\}$. Hence
\bes 
F_x(a,a)(1-a+a^2)=F_y(a,a)(1-a).
\ees
Since $gcd(1-a+a^2,1-a)=1$ it follows from the last equality that
\bes
 F_x(a,a)=r(a)(1-a), \ F_y(a,a)=r(a)(1-a+a^2), \ \text{where} \  r(a) \in \mathbb{C}[a]. 
\ees
 Since $F_x(a,a)\in \mathbb{C}\left[a^2,a^3\right]$ and $(1-a)\notin \mathbb{C}\left[a^2,a^3\right]$ it follows that $r(a) \notin \mathbb{C}^*$.
Consequently, 
\bes
 F=r(A)((1-A)x+(1-A+A^2)y)\pmod{\left\langle a-b\right\rangle}, \ \text{where} \  r(A) \notin \mathbb{C}^*. \Box
\ees

For any automorphism $\psi$ of $M$ we define an ideal
$$ I_{\psi}=\left\langle \psi(h)-h| h\in M\right\rangle. $$
For any $h,m\in M$ we have
\bes
\psi([h,m])-[h,m]=[\psi(h),\psi(m)]-[h,m]\\
=[\psi(h),\psi(m)]-[\psi(h),m]+[\psi(h),m]-[h,m]\\
=[\psi(h),\psi(m)-m]+[\psi(h)-h,m].
\ees
This means that
$$ I_{\psi}=\left\langle \psi(a)-a, \ \psi(b)-b, \ \psi(x)-x, \ \psi(y)-y\right\rangle. $$
Notice that if $M\left\langle a,b,x,y\right\rangle =M\left\langle u ,v,f,g \right\rangle$, then 
$$ I_{\psi}=\left\langle \psi(u)-u, \ \psi(v)-v, \ \psi(f)-f, \ \psi(g)-g\right\rangle. $$

For any automorphisms $\psi$ and $\rho$ of $M$ we also put
 $$ I_{\psi,\rho}=I_{\psi}+I_{\rho}.$$

\begin{lm} \label{lm3} 
We have
\bes
I_{\tau}=\left\langle a-b, t\right\rangle \ \text{and} \ I_{\sigma, \tau}=\left\langle a, b, t\right\rangle, 
\ees
where $$ t=(1+A)x-(1+A+A^2)y.$$
\end{lm}
\Proof It is clear that
$$ I_{\sigma}=\left\langle (w-1)a,  (w^{-1}-1)b,0,0 \right\rangle=\left\langle a,b\right\rangle$$
and
$$I_{\tau}=\left\langle b-a,  a-b, \tau(x)-x, \tau(y)-y \right\rangle=\left\langle a-b, \tau(x)-x, \tau(y)-y \right\rangle.$$
In $M/\left\langle a-b\right\rangle$,  we have
\bes
 \tau(x)-x=-(1+A^3)x+(1+A^2+A^4)y=-(1-A+A^2)t
\ees
and
\bes
\tau(y)-y=(1-A^2)x+(-1+A^3)y=(1-A)t,
\ees
where 
\bes
t=(1+A)x-(1+A+A^2)y.
\ees 
Since $(1-A+A^2)t+A(1-A)t=t $ it follows that $I_{\tau}=\left\langle a-b, t \right\rangle$ and $I_{\sigma,\tau}=I_{\sigma}+I_{\tau}=\left\langle a, b, t\right\rangle.\Box$

\begin{lm} \label{lm4} 
Let $M\left\langle p,q\right\rangle$ be the free metabelian Lie algebra over $\mathbb{C}$ freely generated by $p,q$. Suppose $u$ and $f$ are free  generators of $M\left\langle p,q\right\rangle$ such that 
\bes
u=\alpha p+ P(\mathrm{ad}_{p},\mathrm{ad}_{q})[p,q],\\
f=\gamma q+ Q(\mathrm{ad}_{p},\mathrm{ad}_{q})[p,q],
\ees
where $\alpha,\gamma\in \mathbb{C}^*$, $P(\mathrm{ad}_{p},\mathrm{ad}_{q}), Q(\mathrm{ad}_{p},\mathrm{ad}_{q})\in \mathbb{C}[\mathrm{ad}_{p},\mathrm{ad}_{q}]$. Then
\bes
f=\gamma q+ T(\mathrm{ad}_{p},\mathrm{ad}_{q})[q,[p,q]] \ \  \text{for some} \ \ T(\mathrm{ad}_{p},\mathrm{ad}_{q})\in \mathbb{C}[\mathrm{ad}_{p},\mathrm{ad}_{q}].
\ees
\end{lm}
\Proof Since $u$ and $f$ are free generators of
$M\left\langle p,q\right\rangle$ the endomorphism  
$$\psi=(u,f)=(\alpha p+ P(\mathrm{ad}_{p},\mathrm{ad}_{q})[p,q], \gamma q+ Q(\mathrm{ad}_{p},\mathrm{ad}_{q})[p,q] )$$
 is an automorphism of $M\left\langle p,q\right\rangle$. Hence its Jacobian matrix
\bes
J(\psi)=\begin{bmatrix}
		\alpha-P(p,q)q	 & P(p,q)p \\
		-Q(p,q)q & \gamma+Q(p,q)p
\end{bmatrix},
\ees
$P(p,q), Q(p,q)\in \mathbb{C}[p,q]$,  is invertible over $\mathbb{C}[p,q]$. This means that 
\bes
\det J(\psi)=\alpha\gamma+\alpha Q(p,q)p-\gamma P(p,q)q\in \mathbb{C}^*.
\ees
Therefore, $$ \alpha Q(p,q)p=\gamma P(p,q)q.$$
Setting $q=0$, we obtain $\alpha pQ(p,0)=0$. Since $\alpha\neq 0$ it follows that $Q(p,0)=0$.
Thus $Q(p,q)=T(p,q)q$ for some
$T(p,q)\in\mathbb{C}[p,q]$.
Consequently,
\bes
Q(\mathrm{ad}_{p},\mathrm{ad}_{q})[p,q]
=T(\mathrm{ad}_{p},\mathrm{ad}_{q})[q,[p,q]]. \Box
\ees

\begin{theor} 
The action of $S_3$ on the free metabelian Lie algebra $M$ is not linearizable. 
\end{theor}
\Proof 
Since $\tau(a-b)=b-a$ it follows that $\tau(\left\langle a-b\right\rangle)\subseteq \left\langle a-b\right\rangle$. Then $\tau$ induces an automorphism $\tau'$ of $M/\left\langle a-b\right\rangle\cong M\left\langle a,x,y\right\rangle$ given by
\bes 
\tau'=(a,-A^3x+(1+A^2+A^4)y, (1-A^2)x+A^3y).
\ees

Consider the endomorphism $\varphi=(a,s,t)$ of $M\left\langle a,x,y\right\rangle$, where 
\bes
s=(1-A)x+(1-A+A^2)y  \ \text{and} \ t=(1+A)x-(1+A+A^2)y.
\ees
 The endomorphism $\varphi$ is an automorphism since  
\bes
\det(J(\varphi))=\begin{vmatrix} 1 & 0 & 0 \\ * & 1-a  & 1-a+a^2  \\ * & 1+a &  -(1+a+a^2) \end{vmatrix}=-2.
\ees
	A direct computation shows that $\tau'(a)=a$, $\tau'(s)=s$, $\tau'(t)=-t$ and
\bes
\varphi^{-1}\tau'\varphi=(a,x,-y).
\ees
Consequently, $\tau'$ is linearizable and $M/\left\langle a-b\right\rangle\cong M\left\langle a,s,t\right\rangle$.

Consider the ideals $I_{\tau}$,  $I_{\sigma, \tau}$. By Lemma \ref{lm3},
\bes
I_{\tau}=\left\langle a-b,t\right\rangle \text{ and }  I_{\sigma, \tau}=\left\langle a,b,t\right\rangle .
\ees
It is clear that $I_{\sigma, \tau}\supset I_{\tau}$. Let
$$\overline{M}=M/\left\langle a-b,t\right\rangle=M\left\langle \bar{a}, \bar{s}\right\rangle,$$ 
where  $\overline{h}$ is the image  of $h\in M$ in $\overline{M}$. Then $\overline{I_{\sigma, \tau}}=\left\langle \bar{a}\right\rangle$ in $\overline{M}$. 

The automorphisms $\sigma$ and $\tau$ induce linear
automorphisms
\bes
\tilde{\sigma}=(w\tilde{a},w^{-1}\tilde{b}, \tilde{x},\tilde{y}) \quad \text{and} \quad
\tilde{\tau}=(\tilde{b},\tilde{a},\tilde{y},\tilde{x})
\ees
of the four-dimensional vector space $M/[M,M]$,
where $\tilde{d}$ denotes the image of $d\in M$
in $M/[M,M]$.
As an $S_3$-module,
$\mathbb{C}\tilde{a}\oplus\mathbb{C}\tilde{b}$
is isomorphic to the standard two-dimensional
representation $V$ of $S_3$.
Set
$$
\tilde{h}=\tilde{x}+\tilde{y}
\quad \text{and} \quad
\tilde{m}=\tilde{x}-\tilde{y}.
$$
Since $\tilde{\sigma}(\tilde{h})=\tilde{h}$, $\tilde{\tau}(\tilde{h})=\tilde{h}$,
$\tilde{\sigma}(\tilde{m})=\tilde{m}$, and $\tilde{\tau}(\tilde{m})=-\tilde{m}$,
it follows that $\mathbb{C}\tilde{h}\cong\mathbb{C}$
and $\mathbb{C}\tilde{m}\cong
\mathbb{C}_{\mathrm{sign}}$,
where $\mathbb{C}$ and $\mathbb{C}_{\mathrm{sign}}$
are the one-dimensional trivial and sign representations of $S_3$, respectively. Consequently, we obtain the decomposition
$$
M/[M,M]\cong V\oplus\mathbb{C}
\oplus\mathbb{C}_{\mathrm{sign}}
$$
of the $S_3$-module $M/[M,M]$ into irreducible
submodules.

Suppose that the action of $S_3$ on $M$ is linearizable.
Then there exists a system of free generators
$u,v,f,g$ with respect to which this action is linear. 
Since $u,v,f,g$ are free generators of $M$, their images $\tilde{u},\tilde{v},\tilde{f},\tilde{g}$ form a basis of $M/[M,M]$. Moreover, the restriction of the natural homomorphism $\rho:M\rightarrow M/[M,M]$ to the $S_3$-invariant subspace
$$ W=\mathbb{C}u\oplus\mathbb{C}v\oplus\mathbb{C}f\oplus\mathbb{C}g$$
is $S_3$-equivariant. Hence $\left.\rho\right|_W$
is an isomorphism of $S_3$-modules. Therefore,
$$ W\cong V\oplus\mathbb{C}\oplus \mathbb{C}_{\mathrm{sign}}.$$  Thus, we may choose free generators $u,v,f,g$ of $M$ such that
$$M=M\left\langle a,b,x,y\right\rangle =M\left\langle u,v,f,g\right\rangle$$
and
\bes
\sigma(u)=wu,  \ \sigma(v)=w^{-1}v, \ \sigma(f)=f, \ \sigma(g)=g, \\
\tau(u)=v, \ \tau(v)=u, \  \tau(f)=f,  \ \tau(g)=-g.
\ees
It follows that
$$ I_{\sigma}=\left\langle (w-1)u,  (w^{-1}-1)v, 0, 0 \right\rangle=\left\langle u,v\right\rangle$$
and
$$I_{\tau}=\left\langle u-v,  v-u, 0, 2g \right\rangle=\left\langle u-v, g \right\rangle.$$
Consequently,
\bes
I_{\tau}=\left\langle a-b,t\right\rangle=\left\langle u-v,g\right\rangle  \text{ and } I_{\sigma, \tau}=\left\langle a,b,t\right\rangle=\left\langle u,v,g\right\rangle. 
\ees 
It is clear that $\overline{I_{\sigma, \tau}}=\left\langle\bar{ u}\right\rangle$ in $\overline{M}$. Therefore,
\bee\label{fm1}
\overline{M}=M\left\langle \bar{a},\bar{s}\right\rangle=M\left\langle \bar{u},\bar{f}\right\rangle.
\eee

Denote by $Lin(m)$ the linear part of $m\in M=M\left\langle a,b,x,y\right\rangle$. 

Let $Lin(f)=\alpha a+\beta b+\gamma x+\delta y$. Since $\sigma(f)=\tau(f)=f$ it follows that 
\bes
Lin(\sigma(f))=\alpha w a+\beta w^{-1}b+\gamma x+\delta y=Lin(f),\\
Lin(\tau(f))=\alpha b+\beta a+\gamma y+\delta x=Lin(f).
\ees
Hence $\alpha=\beta=0$ and $\gamma=\delta$. Consequently,
\bes
Lin(f)=\gamma(x+y)=\gamma Lin(s), \ \ \gamma\in \mathbb{C}^*.
\ees
It is easy to show that $Lin(u)=\alpha'a$, $\alpha'\in \mathbb{C}^*$, since $\sigma(u)=wu$. Therefore, by \eqref{fm1}, 
\bes
\overline{u}=\alpha' \overline{a}+ P(\mathrm{ad}_{\overline{a}},\mathrm{ad}_{\overline{s}})[\overline{a},\overline{s}],\\
\overline{f}=\gamma \overline{s}+ Q(\mathrm{ad}_{\overline{a}},\mathrm{ad}_{\overline{s}})[\overline{a},\overline{s}],
\ees
where  $P(\mathrm{ad}_{\overline{a}},\mathrm{ad}_{\overline{s}}), Q(\mathrm{ad}_{\overline{a}},\mathrm{ad}_{\overline{s}})\in \mathbb{C}[\mathrm{ad}_{\overline{a}},\mathrm{ad}_{\overline{s}}]$. By Lemma \ref{lm4}, we have
\bes
\overline{f}=\gamma \overline{s}+ T(\mathrm{ad}_{\overline{a}},\mathrm{ad}_{\overline{s}})[\overline{s},[\overline{a},\overline{s}]] \ \ \text{for some} \ \  T(\mathrm{ad}_{\overline{a}},\mathrm{ad}_{\overline{s}})\in \mathbb{C}[\mathrm{ad}_{\overline{a}},\mathrm{ad}_{\overline{s}}].
\ees
Lifting back to $M$, we have
$$ f=\gamma s+T(\mathrm{ad}_{a},\mathrm{ad}_{s})[s,[a,s]]+h+d$$
for some $ h\in\left\langle t\right\rangle$ and $d\in\left\langle a-b\right\rangle$. Since $s$ and $t$ are homogeneous of degree $1$, and $a$ and $b$ are homogeneous of degree $0$  in $x$ and $y$, the  homogeneous component of $f$ of degree $1$ in $x$ and $y$ has the form
$$ f_1=\gamma s+h_1+d_1, \ \text{ where } h_1\in\left\langle t\right\rangle, \  d_1\in\left\langle a-b\right\rangle.$$
Consequently,
\bee \label{f4}
f_1(a,a,x,y)=\gamma s+h_1(a,a,x,y)= \gamma s+H(A)t, \ \ \text{where} \ \ H(A)\in \mathbb{C}[A].
\eee
Since  $f_1(a,a,x,y)$ is invariant under $\tau'$ and $\tau'(a)=a$, $\tau'(s)=s$, and $\tau'(t)=-t$ it follows from \eqref{f4} that
\bes
\tau'(f_1(a,a,x,y))=f_1(a,a,x,y)=\gamma s-H(A)t. 
\ees
Comparing this with \eqref{f4}, we obtain $H(A)t=0$ and hence $ f_1(a,a,x,y)= \gamma s.$   
This contradicts Lemma \ref{lm2}. $\Box$

\section{An Action of $S_3$ on Metabelian Associative Algebras}

\hspace*{\parindent}

An associative algebra is called {\em metabelian} if it satisfies the identity
\bes
[x,y][z,t]=0.
\ees

Let $A$ be  the free metabelian associative algebra over $K$ with  free generators $x_1,\ldots,x_n$.  
It is well known \cite{UU95}  that the elements of the form 
\bes
x_{i_1}\ldots x_{i_m} [x_{j_k},\ldots, [x_{j_2},[x_{j_1},x_j]] \ldots], 
\ees
where $i_1\leq \ldots \leq i_m$, $j>j_1\leq \ldots \leq j_k$, $m\geq 0$, and $k\geq 1$, is a linear basis of the commutator ideal of $A$. 

Fox derivatives and an analogue of the Jacobian Conjecture for free metabelian associative algebras are studied in \cite{UU95}.

Let $MA=MA\left\langle a,b,x,y\right\rangle$ be the free metabelian associative algebra over the field of complex numbers $\mathbb{C}$ with free generators $a,b,x,y$.  This time denote the operators of the left multiplication $L_a$ and $L_b$ by $A$ and $B$, respectively. Let $C$ be the algebra generated by $A$ and $B$. Obviously, $C$ is isomorphic to the subalgebra of $MA$ generated by $a,b$. 

Let $J$ be the commutator ideal of $C$. Clearly, $J^2 =0$.

Set again
\bes
\sigma=(wa,w^{-1}b,x,y), \ \text{where} \ w=e^{2\pi i/3},
\ees 
and 
\bes
\tau=(
\left[\begin{array}{cc}
	b  \\
	a  \\
\end{array}\right],
S
\left[\begin{array}{cc}
	x  \\
	y  \\
\end{array}\right]),
\ees
 where
 \bes
 S=
 \begin{bmatrix}
 	-B^3+U & 1+AB+A^2B^2+V\\
 	1-AB+W & A^3+Z
 \end{bmatrix}, 
 \ees
$U,V,W,Z$ are arbitrary elements of $J$. Observe that 
\bes
J(\tau)=\begin{bmatrix}
	\begin{matrix}
		0 & 1\\
		1 & 0
	\end{matrix} & \begin{matrix}
		0 & 0 \\
		0 & 0
	\end{matrix} \\
	\begin{matrix}
		* & * \\
		* & *
	\end{matrix} & N
\end{bmatrix}
\ees
has the same form as $J(\tau)$ in the polynomial algebra. Thus $J(\tau)$ is invertible and $\tau$ is an automorphism for any choice of $U,V,W,Z\in J$. 

For any $F=F(A,B)\in C$ set $F^\ast = F(B,A)$. Note that the map $*$ defines an automorphism of $C$. 

\begin{lm}\label{ma1} 
If 
\bes
\begin{bmatrix}
	U & V\\
	W & Z
\end{bmatrix}=-\frac{1}{2}
\begin{bmatrix}
	-B^3 & 1+AB+A^2B^2\\
	1-AB & A^3
\end{bmatrix}
\begin{bmatrix}
	E_{11} & E_{12} \\
	E_{21} & E_{22}
\end{bmatrix},
\ees
where
\bes
E_{11}=[B,A]+B[B,A^2]+[A^3,B^2]B,\\
E_{12}=[B,A^4]+[B^2,A^5],\\
E_{21}=B[A,B^2]B,\\
E_{22}=[A,B]+[A^2,B]B+B[B^2,A^3],
\ees
then $\tau\sigma=\sigma^2\tau$ and $\tau^2=1$.
\end{lm}

\Proof 
Let
\bes
S_0=
\begin{bmatrix}
	-B^3 & 1+AB+A^2B^2\\
	1-AB & A^3
\end{bmatrix} \
\text {and} \
F=
\begin{bmatrix}
	U& V\\
	W& Z
\end{bmatrix}.
\ees
Then $S=S_0+F$ and $F^*F=0$ since $U,V,W,Z\in J$.

Since
\bes
\tau^2=
(\begin{bmatrix}a\\ b\end{bmatrix}, \tau(S\begin{bmatrix}x\\y\end{bmatrix}))
=(\begin{bmatrix}a\\ b\end{bmatrix},S^\ast S\begin{bmatrix}x\\y\end{bmatrix})
\ees
it follows that the equality $\tau^2=1$ holds if 
\bee \label{f5}
S^\ast S=(S^*_0+F^*)(S_0+F)=S^*_0S_0+F^*S_0+S^*_0F=I_2,
\eee
 where  $I_2 $ is the $2 \times 2$ identity matrix.

Let $E=S^*_0S_0-I_2$.  Then
\bes
E=
\begin{bmatrix}
[B,A]+B[B,A^2]+[A^3,B^2]B
&
[B,A^4]+[B^2,A^5]
\\[4pt]
B[A,B^2]B
&
[A,B]+[A^2,B]B+B[B^2,A^3]
\end{bmatrix}
\ees
and, by \eqref{f5},
\bee \label{f6}
F^*S_0+S^*_0F=-E.
\eee

Notice that 
\bes
E^\ast S_0=(S_0S_0^\ast-I_2)S_0=S_0(S_0^\ast S_0-I_2)=S_0E
\ees
and $E^2=0$ since all elements of $E$ belong to $J$. 
 Using these we have
\bes
-\frac{1}{2}(S_0^\ast E^\ast S_0+S_0^\ast S_0 E)=- S_0^\ast S_0 E=-(E+I_2)E=-E.
\ees
Then 
\bes
F=-\frac{1}{2} S_0 E
\ees
is a solution of equation \eqref{f6}.

We have 
\bes
\tau\sigma=(\begin{bmatrix}wb\\w^{-1}a\end{bmatrix},
S_0\begin{bmatrix}x\\y\end{bmatrix}+ F \begin{bmatrix}x\\ y \end{bmatrix})
\ees
and 
\bes
\sigma^2\tau=(\begin{bmatrix}wb\\w^{-1}a\end{bmatrix},
S_0 \begin{bmatrix}x\\y\end{bmatrix}
+\b(F) \begin{bmatrix} x\\ y \end{bmatrix}), 
\ees
where $\b$ is the automorphism of $C$ induced by $\sigma^2=(w^{-1}a, wb, x, y)$.
It is not difficult to see that $\b(F)=F$. Consequently, $\tau \sigma=\sigma^2 \tau$. $\Box$

Let $\tau$ denote the automorphism of $MA$ determined by $U,V,W,Z$, which satisfy the assumptions of Lemma \ref{ma1}:
\bes
\tau=(\begin{bmatrix}b\\a\end{bmatrix},
\begin{bmatrix}
	-B^3 & 1+AB+A^2B^2\\
	1-AB & A^3
\end{bmatrix}
\begin{bmatrix}x\\y\end{bmatrix}-\frac{1}{2} 
\begin{bmatrix}
	-B^3 & 1+AB+A^2B^2\\
	1-AB & A^3
\end{bmatrix}
\begin{bmatrix}
	E_{11} & E_{12} \\
	E_{21} & E_{22}
\end{bmatrix}
\begin{bmatrix}x\\ y \end{bmatrix}),
\ees 
where
\bes
E_{11}=[B,A]+B[B,A^2]+[A^3,B^2]B,\\
E_{12}=[B,A^4]+[B^2,A^5],\\
E_{21}=B[A,B^2]B,\\
E_{22}=[A,B]+[A^2,B]B+B[B^2,A^3].
\ees

\begin{co}
The automorphisms $\sigma$ and $\tau$ of $MA$ satisfy the relations (\ref{f1}). 
\end{co}

\begin{theor}\label{ma2} The action of $S_3$ on the free metabelian associative algebra $MA$ is not linearizable. 
\end{theor}
\Proof Observe that the action of $S_3$ on $MA$ induces the action of $S_3$ on the polynomial algebra $\mathbb{C}[a,b,x,y]$ described in Section 2. Since the latter is not linearizable, it follows that the action of $S_3$ on $MA$ is also not linearizable. 
$\Box$

\section*{Acknowledgments}

This work was supported by the grant AP23487886 of the Ministry of Science and Higher Education of the Republic of Kazakhstan.

We used AI tools to improve the grammar of the manuscript. In the initial version of the paper, we also used AI to solve certain challenging systems of linear equations. In the present version, however, we have simplified the calculations  without using AI.

\end{document}